\documentclass[a4paper,12pt]{article}
\input{preamble}

\newcommand{\jpTitle}{Support in the domain of numerical algorithms in the education of future secondary school mathematics teachers}
\newcommand{\jpAuthors}{J. Dan\v{e}k}
\newcommand{\jpKeywords}{Digital educational resources; mathematics education; STEM; numerical methods;  MATLAB}

\author{Josef Dan\v{e}k\thanks{\email{danek@kma.zcu.cz}}}
\affil{Department of Mathematics, Faculty of Applied Sciences, \authorcr University of West Bohemia in Pilsen, Univerzitn\'{\i} 2732/8, 301 00 Plze\v{n}, Czech Republic}

\title{\textcolor{Navy}{\textsc{\jpTitle}}}

\begin{document}

\maketitle

\begin{abstract}
This article presents web-based applets designed primarily to educate future mathematics teachers in the field of numerical algorithms and, secondarily, to serve high school students themselves. The applets allow users to explore the behaviour of basic numerical algorithms as applied to fundamental mathematical topics taught at high schools. The advantage is that no software is required, and students can intuitively test how a given method works by adjusting specific problem parameters, not only by viewing a list of approximate solutions, but especially through graphical output. Furthermore, they have access to multiple different methods for comparison.
\end{abstract}


\section{Introduction}\label{sec:introduction}
Teaching and learning focused on the transfer of mathematical skills and knowledge presents a challenge. Large classes, individual needs and the varying programmes of individual pupils, together with a generally negative attitude towards mathematics, make this challenge even harder. Therefore, it is necessary to engage in research on teaching methods and technologies that focus on the transmission of mathematical knowledge and skills \cite{Ari2020, Cullen2020, Koyunkaya2024}.

Mathematics taught at secondary school typically includes topics such as finding roots of polynomials, computing eigenvalues of a matrix, calculating derivatives of functions, and evaluating definite integrals. These topics are practiced using examples where the solution can be determined exactly by a simple direct method. However, when it comes to finding, for example, the roots of a general polynomial of degree five or higher, an exact solution is generally not possible, and one must resort to approximate numerical methods. Such methods are not usually part of the secondary school curriculum, but we believe that secondary school teachers should be able to explain to students the fundamental ideas behind numerical algorithms, that these problems could not be solved without them.

This article is based on many years of experience teaching mathematics (mathematical analysis, numerical mathematics, linear algebra) especially at the technical faculties of the University of West Bohemia in Pilsen (students aged 19-24), which gave me information about the state of mathematical knowledge of students entering the first year. In the last six years I have also had the opportunity to teach mathematics at a technical high school (students aged 15-19) and thus I have tried to prepare these pupils for future university studies. Thanks to the experience gained, I can now specifically concentrate on students of the Faculty of Education studying to become mathematics teachers for secondary schools and attempt to pass on as much knowledge and information as possible.

Important factors that currently influence teaching have been and are the COVID-19 pandemic and new digital technologies, especially the use of generative AI. Greater emphasis is placed on the potential to help students study independently at their own pace. Digital technologies have been used in teaching before, but the present situation forces us to think more about the resources we use and leads to an attempt to expand the range of teaching materials \cite{Liu2021,Weinhandl2022},
for example by using pre-recorded videos and communication tools (Zoom, Google Meet) 
\cite{Chen2021,Cevikbas2025} 
and efficiently sharing knowledge with students \cite{Zochbauer2021}. 

In previous articles I have focused on the use of digital technologies in the teaching of numerical
 mathematics at university \cite{Danek2004,Danek2006}. 
 Over the past two decades, I have had the opportunity to directly compare mathematics instruction in high schools with the training of future mathematics teachers at the University of West Bohemia.
 This article shows the possibilities of using mathematical software in the training of future mathematics teachers in the area of basic numerical methods and of course subsequently in their own practice for teaching secondary school students. All the applets created and presented here have been used in teaching, and according to my experience the students appreciated the clear illustration of basic numerical mathematics algorithms on simple examples with the possibility of use outside contact teaching and with their own pace of practice in the form of independent work.

The aim of the course ``Seminar on Numerical Methods for Teachers'' which is the focus of this article, is to provide knowledge about algorithms of approximate computational methods and to teach how to convey this knowledge to secondary school students. However, the intention is not merely to train students to use tools; rather, it is to help them understand the underlying principles and to make the processes of teaching and learning more efficient and of higher quality \cite{Maddux2011}. 

The applets were created using MATLAB, a specialised software environment for numerical computation, data analysis, simulation, and programming, which is widely used in science, engineering, mathematics, and technology. MATLAB was developed by Cleve Moler in the late 1970s; the commercial version was released in 1984, when the company MathWorks (https://mathworks.com/) was founded. Today, MATLAB is among the most widely used tools for technical computing and simulation worldwide. Another tool for creating digital teaching materials is the dynamic software GeoGebra, which combines the features of dynamic geometry software and computer algebra systems. One example of its use in teaching is the geometric design of building roofs \cite{Tomiczkova2025}. 

The use of digital tools in mathematics education can significantly improve students' ability to design and solve real-world mathematical problems \cite{Greefrath2018}. 
Their successful use by future teachers depends primarily on their digital competences, the availability 
of technologies, and the support provided by educational institutions 
\cite{Ishartono2025,Godfrey2008,Yu2026}. 
It must be recognised that technology alone is not sufficient; what is crucial is its pedagogically meaningful integration into teaching \cite{Greefrath2011,Weintrop2016}.  

\section{Methodology}\label{sec:methodology}
The standard secondary school mathematics curriculum includes manipulation of expressions, sets, linear and quadratic equations, inequalities and their systems, trigonometry, goniometry, solid geometry, analytic geometry, elementary functions, sequences, series, combinatorics, probability, statistics, complex numbers, matrices, limits of functions, derivatives, and definite and indefinite integrals. For the presentation and practice of this content, examples with ``nice solutions'' are typically used, for instance equations with integer solutions that can easily be found by a direct method. This is of course desirable for an initial understanding of the topic. However, in the Czech Republic, in addition to general secondary schools, there are also schools more oriented towards the future profession of the graduate, such as secondary technical schools specialising in mechanical engineering, civil engineering, or electrical engineering. Especially at such technically oriented schools, I consider it important that students at least briefly encounter the basic ideas of numerical algorithms by which these problems are solved in practice. If they use specialised software in vocational subjects, they should at least have a general idea of how this software works, rather than treating it merely as a ``black box''.

Secondary school mathematics teachers' approaches to using digital technologies at different stages 
of the mathematical modelling process are examined in a study by 
\cite{Kandemir2025}. 
The results show that teachers used the internet, spreadsheets, calculators, and mathematical software effectively to solve problems. However, the challenges included excessive reliance on technological outputs, limited critical evaluation and validation of models, and difficulties in transferring mathematical knowledge into the modelling process. These findings highlight the need for explicit support for instruction in teacher education programmes, in order to strengthen teachers' critical engagement with digital tools and enhance their ability to integrate mathematical knowledge into real-world problem solving.

A variety of tools can be used to illustrate numerical algorithms, for example MATLAB. 
The results of the analysis by \cite{Wang2024} 
indicate that the use of MATLAB increases the visual clarity of teaching advanced mathematics and can contribute to improving both teaching quality and student motivation. 
To avoid software installation, this article presents applets that are created in MATLAB but are executable and accessible via a web browser.

For mathematics teachers to be able to teach mathematics effectively in schools, they must possess
 appropriate knowledge and competences \cite{Guerrero-Ortiz2021}.  
 It is necessary to connect mathematical knowledge, didactic skills, and the ability to work with real-world problems. My experience with teaching future teachers shows that, in their mathematical training, emphasis is placed primarily on the secondary-school curriculum itself. Consequently, they lack the broader mathematical perspective that is important for communicating mathematical knowledge. This deficiency is addressed in the course ``Seminar on Numerical Methods for Teachers'' in which the applets described here can be advantageously employed.

\section{Structure of the applets}\label{sec:structureoftheapplets}
The applets presented here offer an effective tool for understanding the principles, properties, and behaviour of basic numerical methods. Through an interactive user interface, they enable users to experiment with different input parameters, immediately observe the results, visualise the convergence of algorithms or the behaviour of errors, and thus support a faster understanding of the characteristics of individual algorithms and their mutual comparison. Owing to their user-friendliness and immediate visual feedback, they enrich the educational process.  
The applets focus on numerical methods for solving, in particular: nonlinear equations, systems of linear algebraic equations, eigenvalues and eigenvectors, function approximation, numerical differentiation, and definite integration.

\begin{figure}[!ht]
\begin{center}
\includegraphics[width=0.98\textwidth]{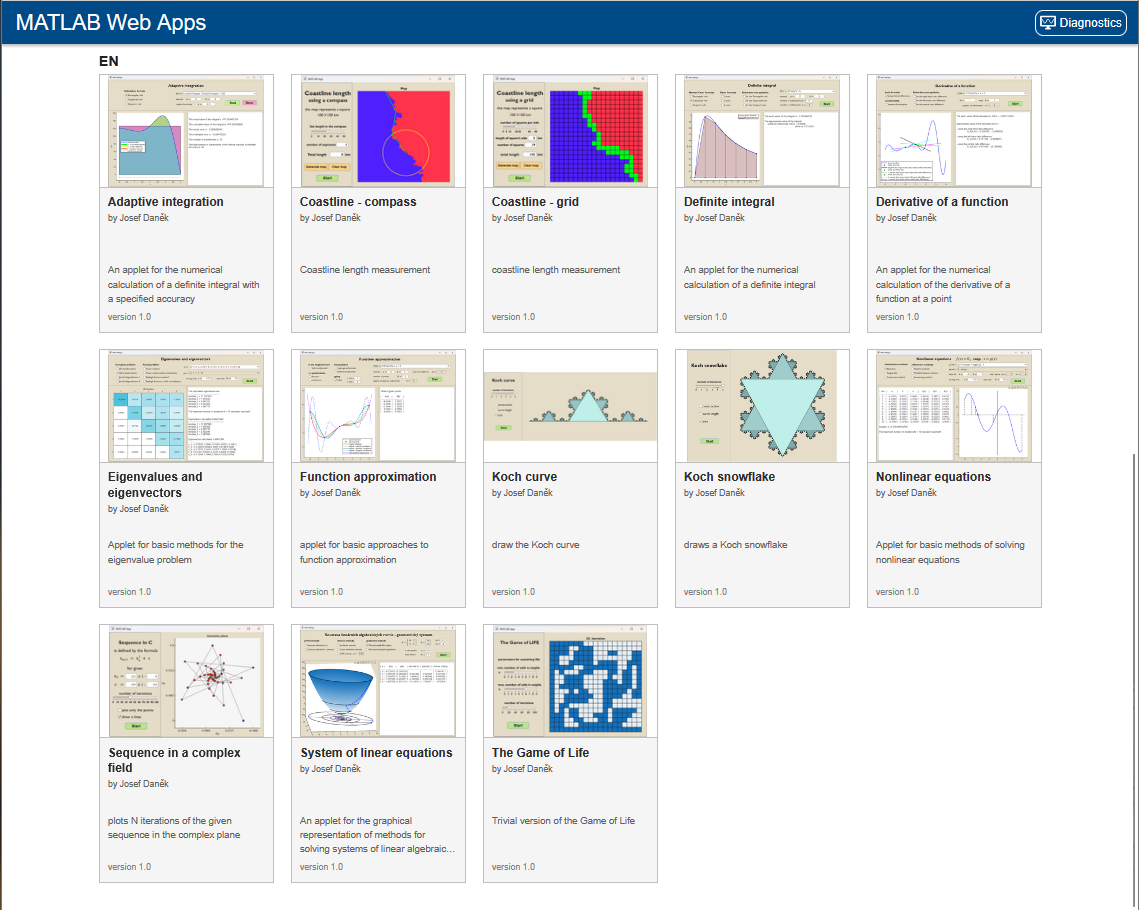}
\end{center}
\caption{List of applets on the MATLAB web server \href{https://matlab.kma.zcu.cz}{https://matlab.kma.zcu.cz}}\label{f:matlabweb}
\end{figure}

All applets (available at \href{https://matlab.kma.zcu.cz}{https://matlab.kma.zcu.cz}) share the following structure: the upper part contains controls for selecting the method, specifying a particular example, and entering individual parameters. The larger lower part is divided into two windows. The first provides a graphical output describing the solution for the chosen problem and method, while the second window is used, for example, to display the individual iterations and the numerical solution of the method.  
These windows are interconnected: for instance, in an iterative method, the approximate solution is plotted step by step in the first window, while the corresponding numerical values are simultaneously listed in the second. The graphical output supports zooming in and out and, in the case of 3D plots, also rotation, which helps users to understand the problem more quickly.

\subsection*{Nonlinear equations}
\begin{figure}[!ht]
\begin{center}
\includegraphics[width=0.98\textwidth]{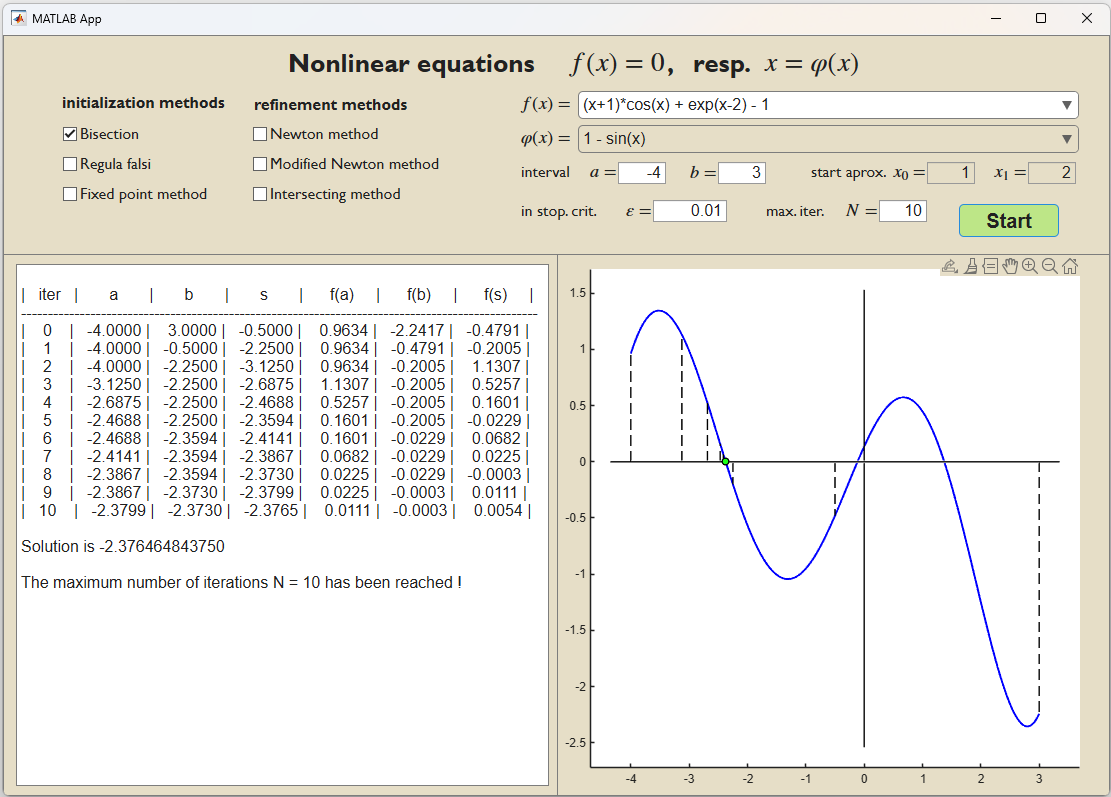}
\end{center}
\caption{Applet for solving nonlinear equations}\label{f:nelinearnirovnice}
\end{figure}

In this example (see Figure \ref{f:nelinearnirovnice}), we solve the equation  
$$
(x+1)\cdot \cos x + \e^{x - 2} -1 = 0
$$
on the interval $[ -4; 3 ]$  using the bisection method. In the right-hand window, the graph of the given function is plotted at each iteration with the current interval highlighted, while in the left-hand window the corresponding iterations of the chosen method are listed step by step. It is evident that there are three distinct solutions on the given interval, so at the beginning the question arises as to which root the method will converge.  
Additional methods are available for solving the problem: the Regula Falsi method, fixed-point iteration, Newton's method, the modified Newton method, and the secant method. By changing the parameters, we can explore the behaviour of each method individually and also compare them with one another (in terms of convergence, number of iterations, and error).

\subsection*{System of linear equations}
\begin{figure}[!ht]
\begin{center}
\includegraphics[width=.98\textwidth]{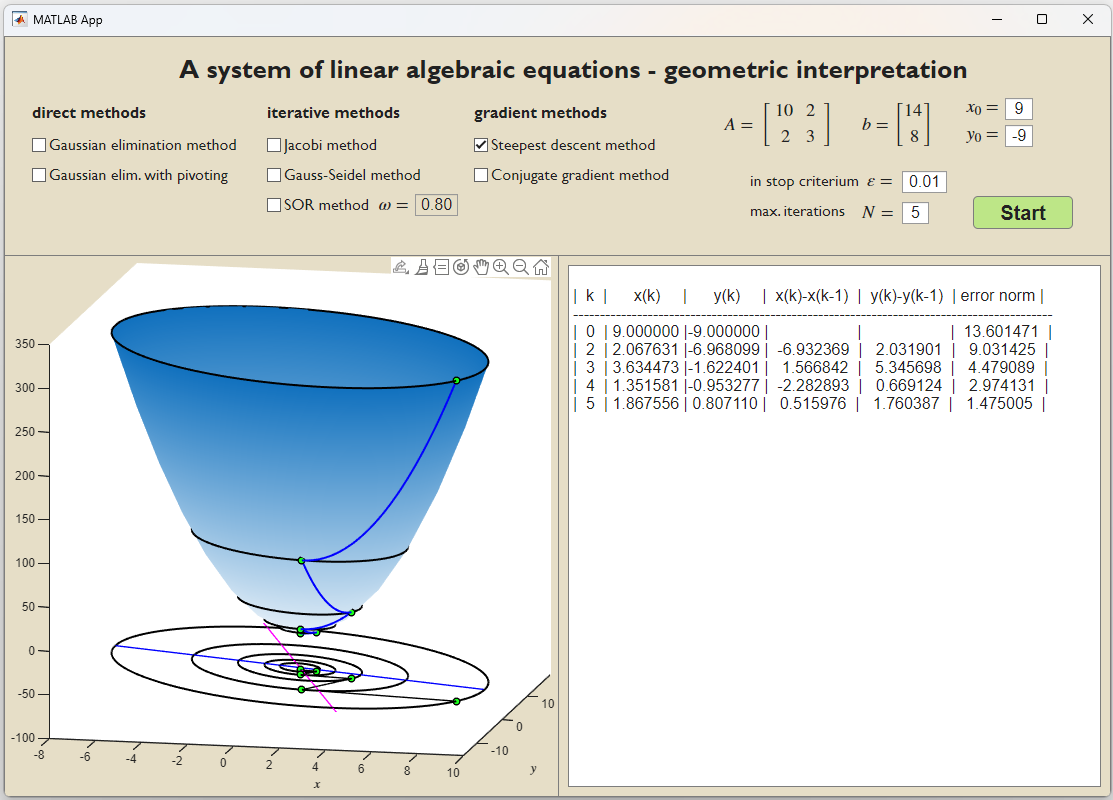}
\end{center}
\caption{Applet for methods of solving systems of linear equations}\label{f:slar}
\end{figure}

In this example (see Figure \ref{f:slar}), we solve the system  
$$
\begin{array}{c} 
10 x +2 y = 14 \\ 
2 x +3 y = 8\\  
\end{array}
$$
using the method of steepest descent. The initial approximation is chosen as 
$[x_0, y_0]^T = [9, -9]^T$, and five iterations are computed and displayed.  
From the graph in the left window we gain a clear insight into the idea of gradient methods, namely the minimisation of the corresponding quadratic function, whose graph is an elliptic paraboloid. The entire iterative process can be observed in 3D and simultaneously in its projection onto the x,y-plane, where the corresponding level curves are plotted and the alternation of orthogonal directions in successive iterations becomes apparent. The numerical values of the iterations are listed step by step in the right-hand window.  
For comparison, it is also possible to use the conjugate gradient method, the iterative methods of Jacobi, Gauss-Seidel, and successive over-relaxation, as well as the direct Gaussian elimination method with and without pivoting.

\subsection*{Eigenvalues and eigenvectors}
\begin{figure}[!ht]
\begin{center}
\includegraphics[width=.98\textwidth]{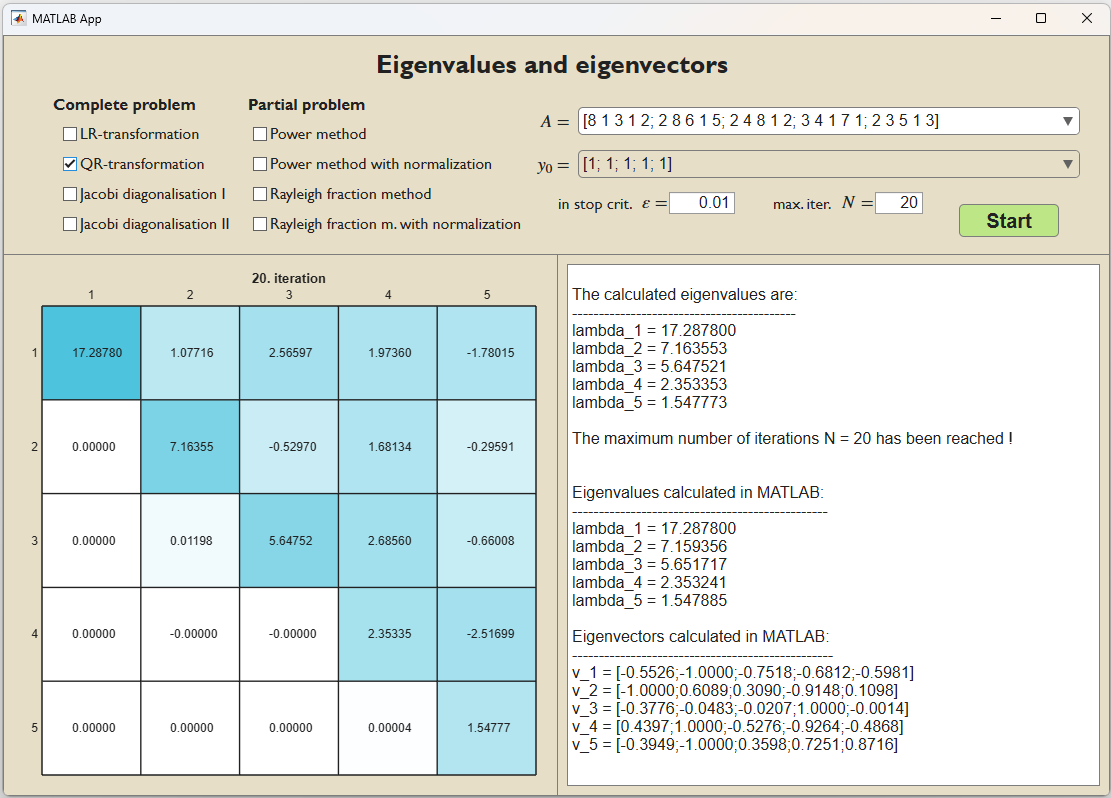}
\end{center}
\caption{Applet for methods of finding eigenvalues and eigenvectors}\label{f:vlastnicisla}
\end{figure}

In the sample example (see Figure \ref{f:vlastnicisla}), all eigenvalues of the matrix 
$$
{\mathbf A} = \begin{bmatrix}
8 & 1 & 3 & 1 & 2 \\ 2 & 8 & 6 & 1 & 5 \\ 2 & 4 & 8 & 1 & 2 \\ 3 & 4 & 1 & 7 & 1 \\ 2 & 3 & 5 & 1 & 3
\end{bmatrix}
$$
are computed using the QR-iteration method, in which a sequence of pairwise similar matrices (i.e. with identical eigenvalues) is generated that converges to an upper triangular matrix with the eigenvalues on its diagonal. In the left window, the individual matrix iterates are displayed, and for clarity each entry is colour-coded with an intensity corresponding to the magnitude of its absolute value. In the right window, the computed eigenvalues are listed, together with the eigenvalues and eigenvectors obtained using the MATLAB function \verb|eig|.  
For the full problem (computing all eigenvalues), the LR-iteration method is also available, and for symmetric matrices additionally two versions of the Jacobi diagonalisation method. For the partial problem, where only the eigenvalue of largest absolute value is sought, the Power Method and the Rayleigh quotient method, both with and without normalisation, are provided.

\subsection*{Function approximation}
\begin{figure}[!ht]
\begin{center}
\includegraphics[width=.98\textwidth]{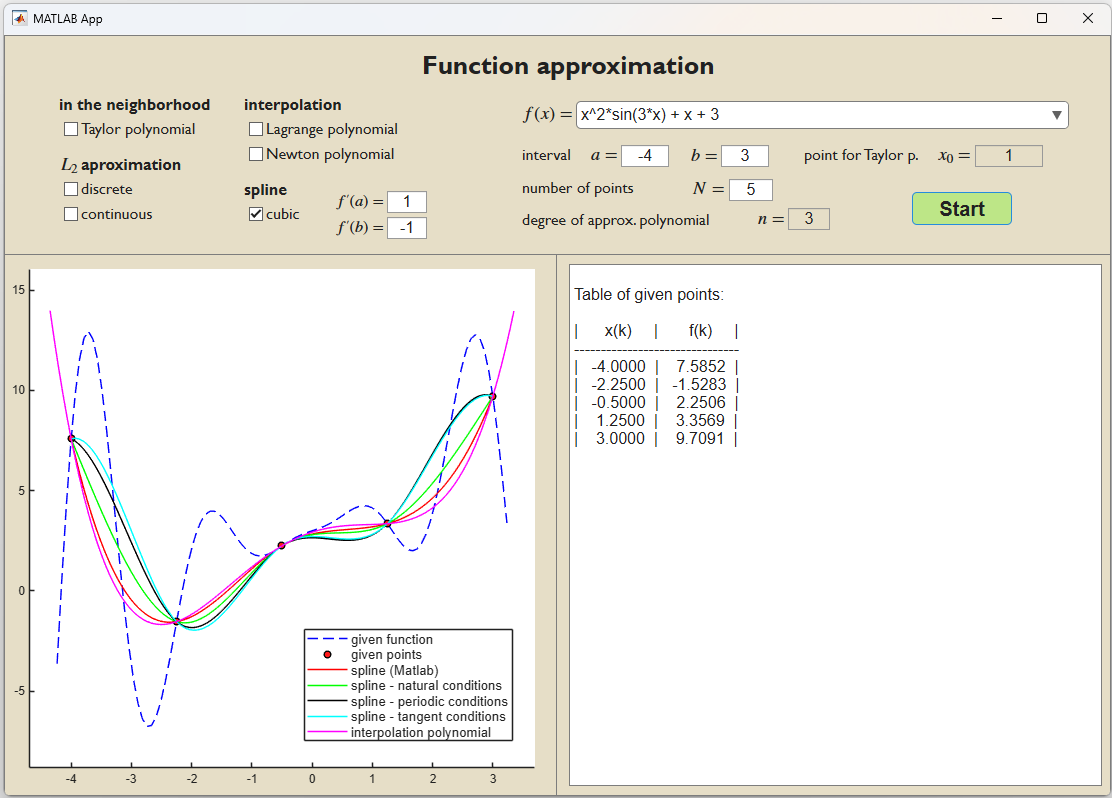}
\end{center}
\caption{Applet for various approaches to function approximation}\label{f:aproximacefunkce}
\end{figure}

In the sample example (see Figure \ref{f:aproximacefunkce}), the function  
$f(x)=x^2\cdot \sin(3x) +x +3$
is given at five points chosen equidistantly on the interval $[ -4; 3 ]$. For the approximation, a cubic spline is used in four variants corresponding to different additional conditions (no explicit conditions - MATLAB supplies them via interpolation, natural boundary conditions, periodic boundary conditions, and clamped boundary conditions), and for comparison the interpolating polynomial is also plotted.  
For interpolation, either the Lagrange or the Newton interpolating polynomial can be used. In addition, the Taylor polynomial can be plotted and printed. The final type of function approximation implemented is the $L_2$-approximation, available in both discrete and continuous variants.

\subsection*{Derivative of a function}
\begin{figure}[!ht]
\begin{center}
\includegraphics[width=.98\textwidth]{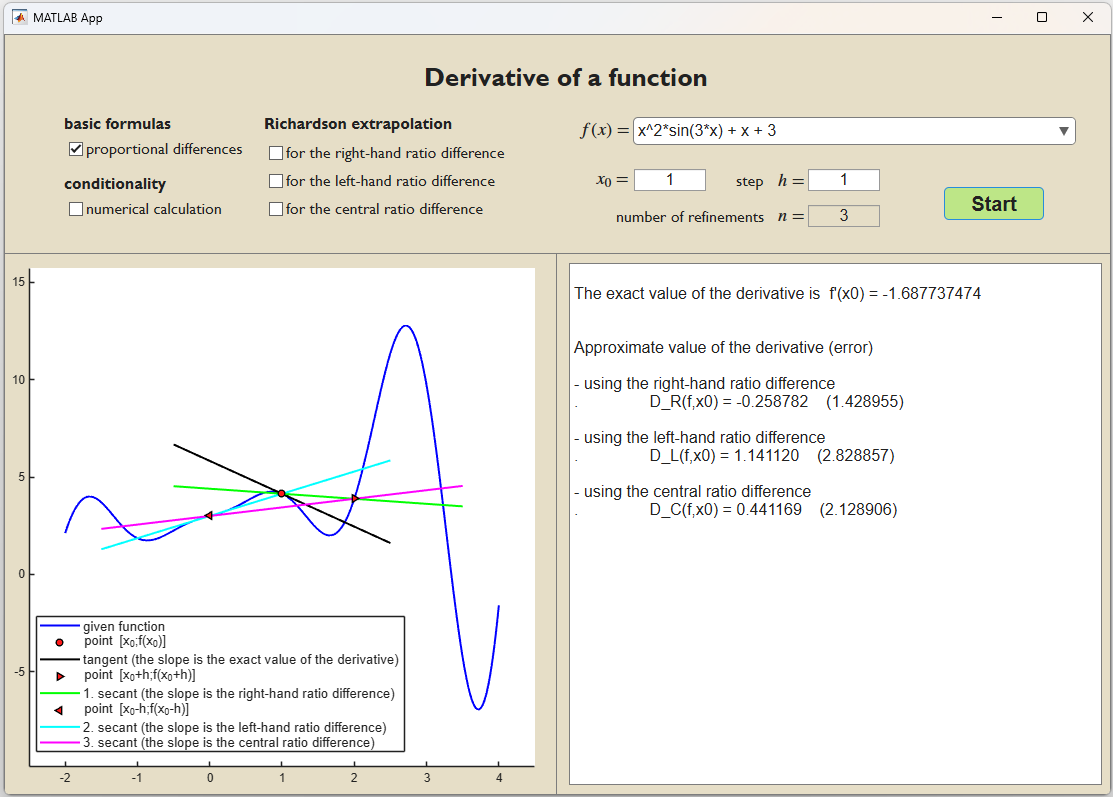}
\end{center}
\caption{Applet for approximately calculating the derivative of a function at a point}\label{f:derivacefunkce}
\end{figure}

In the sample example (see Figure \ref{f:derivacefunkce}), the geometric meaning of the basic formulas for computing the derivative of a function is illustrated, namely the forward, backward, and central finite differences. Specifically, we again consider the function $f(x)=x^2\cdot \sin(3x) +x +3$,
with the derivative evaluated at the point $x_0=1$ and step size $h=1$. In the left window, the graph of the given function is plotted together with the corresponding lines whose slopes represent the derivative approximations. In the right window, the exact value of the derivative is given and the numerical results of the formulas, including the errors, are displayed.  
The applet also makes it possible to illustrate the ill-conditioning of the numerical differentiation problem and to use Richardson extrapolation to improve accuracy.

\subsection*{Definite integral}
\begin{figure}[!ht]
\begin{center}
\includegraphics[width=.98\textwidth]{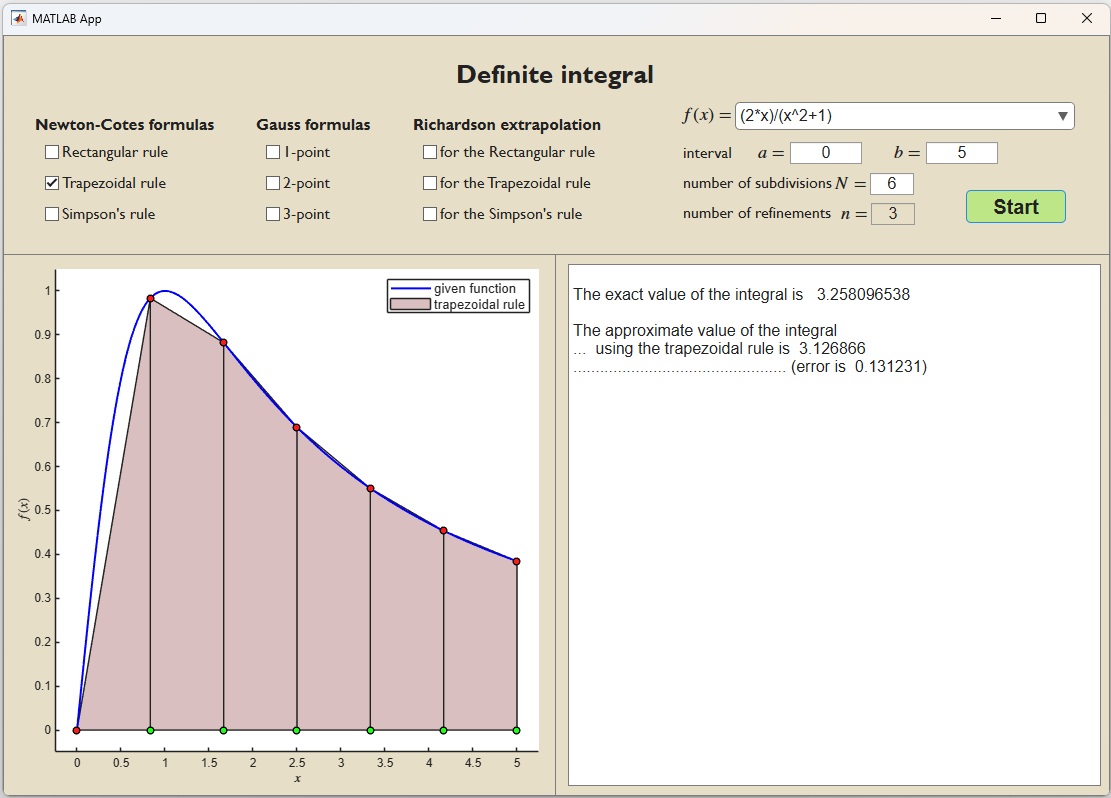}
\end{center}
\caption{Applet for methods of calculating a definite integral}\label{f:urcityintegral}
\end{figure}

Numerical methods for calculating the value of a definite integral are based on its geometric meaning 
(the area between the graph of a given function and the x-axis).
In the example used (see Figure \ref{f:urcityintegral}), we compute an approximate value of the definite integral  
$$\int\limits_0^5 \frac{2x}{x^2+1} \,\d x$$ 
using the composite trapezoidal rule with step size 1. In the left window, the graph of the integrand is plotted together with the individual trapezoids corresponding to the chosen setting. In the right window, the exact value and the approximate value obtained by the selected formula are displayed, including the error.  
The list of methods further includes other Newton-Cotes formulas (the rectangle rule and Simpson's rule) as well as composite Gaussian quadrature formulas. For the Newton-Cotes formulas, Richardson extrapolation can be used to improve accuracy.

\section{Conclusion}\label{sec:conclusion}
University-level education of future secondary school mathematics teachers in the area of basic numerical algorithms is absolutely crucial in today's digital era and in the context of the ever-growing importance of data and computation. Algorithms for solving nonlinear equations, systems of linear equations, computing matrix eigenvalues, function approximation, numerical differentiation, or evaluating definite integrals
form a bridge between abstract mathematical theory and its practical applications. For prospective teachers, it is essential not only to master analytical methods, but also to understand when and why to resort to numerical approaches, how these algorithms work internally, and what their limitations are. 

A solid ground in numerical mathematics enables them to explain complex concepts more effectively and clearly, to demystify the ``black boxes'' of modern calculators and software, and to demonstrate practical applications of mathematics in the real world: from engineering and economics to biology. Such preparation is indispensable for fostering students' critical thinking about computer-based computation, developing their ability to tackle complex problems, and preparing them for further study at university, where numerical methods are encountered extensively. Acquiring these competences is therefore fundamental for educating qualified teachers who can inspire pupils and equip them with skills essential for the 21st century.

A wide range of software tools and applications is available for the presentation and demonstration of numerical methods, differing in complexity, purpose, and user friendliness. When making a selection, it is important to consider the target group (secondary school students), the complexity of the method being demonstrated, and the desired level of interactivity and visualisation. The tools can be roughly divided into several categories:
\begin{itemize}
\item Programming languages and libraries - for deeper understanding and flexibility (Python with NumPy, SciPy and Matplotlib, MATLAB, GNU Octave, Julia)
\item Mathematical software - for advanced symbolic and numerical computation and visualisation (Wolfram Mathematica, Maple, Sage Math)
\item Dynamic geometry software and graphing calculators (GeoGebra, Desmos)  
\item Spreadsheets (Microsoft Excel, Google Sheets, LibreOffice Calc)
\item Online calculators and simulators (Wolfram Alpha, Java applets or interactive web applications)
\end{itemize}
However, some tools are not well suited for classroom use, mainly because of licensing constraints, cost, and especially software complexity. Individual tools also differ in their focus on particular types of problems.

For the materials presented here, MATLAB was used to create web-based applets that can be run directly in a browser, without the need for a local MATLAB installation (which requires a licence). The advantages are easy accessibility, no installation, and the possibility of quickly demonstrating specific algorithms. The materials are useful for both distance and face-to-face teaching and are particularly appreciated by students for independent study. They have been made available on a public website, 
\href{https://matlab.kma.zcu.cz}{https://matlab.kma.zcu.cz}, and can thus be used by students and teachers at other institutions as well.

The materials can be beneficial on two levels. First, they can be employed in university programmes for future secondary school mathematics teachers, who, through these applets, become familiar with the ideas and principles of numerical methods for the topics covered in secondary school. Subsequently, these teachers can use the same tools in their own teaching practice and pass on the acquired knowledge and skills to their secondary school students. In this way, they contribute to improving the quality and effectiveness of mathematics teaching and to a better understanding of numerical methods, supported by direct visualisation of the results.
 
\section*{Disclosure statement} 
No potential conflict of interest was reported by the author.

\section*{Funding}
No funding was received.

\section*{ORCID}
Josef Dan\v{e}k  \href{https://orcid.org/0000-0002-9744-5107}{https://orcid.org/0000-0002-9744-5107}.

\renewcommand*{\bibfont}{\normalsize}

\end{document}